\documentclass[12pt]{article}
\usepackage{amssymb,latexsym,amsmath,eufrak,amsthm,cite,pictex}
\usepackage[active]{srcltx}
\usepackage{hyperref}
\usepackage[inline]{enumitem}
\usepackage[T1]{fontenc}
\usepackage{mathrsfs} 
\allowdisplaybreaks[1]

\def\gH{\EuFrak H}

\def\gM{\EuFrak M}

\def\vphi{\varphi}

\def\l{\lambda}

\def\s{\sigma}
\def\t{\tau}
\def\om{\omega}

\def\Om{\Omega}

\def\z{\zeta}

\newcommand{\D}{{\mathbb D}}
\newcommand{\T}{{\mathbb T}}

\makeatletter
\newcommand{\inlineitem}[1][]{%
\ifnum\enit@type=\tw@
    {\descriptionlabel{#1}}
  \hspace{\labelsep}
\else
  \ifnum\enit@type=\z@
       \refstepcounter{\@listctr}\fi
    \quad\@itemlabel\hspace{\labelsep}
\fi}
\makeatother

\newcommand{\ol}[1]{{\overline{#1}}}

\theoremstyle{plain}
\newtheorem{theorem}{Theorem}
\newtheorem{corollary}{Corollary}
\newtheorem{lemma}{Lemma}

\theoremstyle{definition}
\newtheorem{definition}{Definition}
\theoremstyle{remark}
\newtheorem{remark}{Remark}
\title{On Generators of the Hardy and the Bergman Spaces}
\author{Valentin~V.~Andreev, Miron~B.~Bekker\footnote{Corresponding author}, Joseph~A.~Cima}
\date{}
\begin{document}
\maketitle
\begin{abstract}
A function $\vphi$ which is analytic and bounded in the unit disk $\D$ is called a generator for the Hardy space $H^2(\D)$ or the Bergman space $A^2(\D)$ if polynomials in $\vphi$ are dense in the corresponding space. We characterize generators in terms of $\vphi-$invariant subspaces which are also $z-$invariant and study wandering properties of such subspaces. Density of bounded analytic functions in the $\vphi-$invariant subspaces  is also investigated.
\end{abstract}



\section{Introduction}
\label{}
It is well known that any (closed) subspace of the Hardy space $H^2(\D)$ or the  Bergman space $A^2(\D)$ which is invariant under multiplication by $z$ (the $z-$invariant subspace) is also an invariant subspace of any bounded analytic Toeplitz operator. The converse is not in general true.

In this article we consider analytic Toeplitz operators $M_{\vphi}$ of multiplication by $\vphi$ on the Hardy space $H^2(\D)$ and the Bergman space $A^2(\D)$ with the property that the function which is a constant (e.g. $f(z)\equiv 1$) is a cyclic vector for such operators.  Such a function $\vphi$ we call a generator. We characterize such generators in terms of $M_{\vphi}-$ invariant subspaces which are also $z-$invariant. 

The problem of determining whether a given bounded analytic function in the the unit disk $\D$ is a generator is very hard. Some important results in this direction in the case of Hardy space were obtained in the paper by B.~M.~Solomyak \cite{Solomyak}.

In \cite{Sarason_1} and \cite{Sarason_2} Donald Sarason introduced a notion of a weak$^*$ generator of $H^{\infty}$. We show that each weak$^*$ generator is a generator in our sense.  If $\psi$ is a weak$^*$ generator of $H^{\infty}$ then $\psi-$invariant subspaces have some additional property, see Theorem \ref{H-infty dense} below. 

Weak$^*$ generators often occur in the process of studying wandering property of invariant subspaces of analytic Toeplitz operators. 
In this article we also consider this question for the case when the symbol of the Toeplitz operator is a generator, see Theorem \ref{wandering_property_th}. 

The condition that $H^{\infty}\cap\gM$ is dense in the invariant subspace $\gM$  plays the crucial role in our work.  In  Section 3 we study this condition in the case of invariant subspaces of the Hardy space $H^2(\D)$ and the Bergman space $A^2(\D)$.
\section{Generators}
In what follows $\D$ denotes the unit disk and $\gH$ denotes either Hardy space $H^2(\D)$ or Bergman space $A^2(\D)$.

\begin{definition}\label{generator} A function $\vphi\in H^{\infty}(\D)$ is called  a  generator for $\gH$ if polynomials in $\vphi$ are  dense in $\gH$, i.e.
\begin{equation}\label{clh}
 \overline{l.h.\{\vphi^n:n=0,1,\ldots\}}=\gH
\end{equation}
where {\it l.h.} means linear hull and bar means closure in the norm of $\gH$. 
\end{definition}
Condition \eqref{clh}  implies that a vector ${\mathbf 1}$ (the constant function $1$) is a cyclic vector for the analytic Toeplitz operator $M_{\vphi}$ of multiplication by $\vphi$. 

Since in a locally convex space the weak closure of a subspace coincides with its norm closure, a function $\vphi$ is a generator if polynomials in $\vphi$ are weakly dense in $\gH$. 

Clearly, $\vphi$ is univalent. Indeed, if $\vphi(z_1)=\vphi(z_2)$, $z_1,z_2\in\D$, $z_1\ne z_2$, then 
$P(\vphi(z_1))=P(\vphi(z_2))$ for any polynomial $P$. Consequently, $f(z_1)=f(z_2)$ for any $f\in\overline{l.h.\{\vphi^n:n=0,1,\ldots\}}$  which contradicts \eqref{clh}.

A notion of a generator is closely related to the notion of weak$^*$ generator introduced by D.~Sarason in  \cite{Sarason_1}, \cite{Sarason_2}. A function 
$\psi\in H^{\infty}$ is called a weak$^*$ generator of $H^{\infty}$ if polynomials in $\psi$ are dense in the weak-star topology of $H^{\infty}$. 
\begin{theorem}
 \begin{enumerate}
\item Every weak$^*$ generator is a generator in the sense of Definition \ref{generator}. 
\item There exist a bounded univalent function $\vphi$ which is a generator in the sense of Definition \ref{generator} but is not a weak$^*$ generator.
\end{enumerate}
\end{theorem}
{\bf Proof.} {\it 1.} Recall that the space $H^{\infty}$ is the dual to the quotient space $L^1/H^1_0$, $H^{\infty}=(L^1/H_0^1)^*$.  A local basis of the weak$^*$ topology on $H^{\infty}$ at $h_0\in H^{\infty}$ 
is formed by the following sets
\begin{equation}\label{basis_topology}
\left\{h\in H^{\infty}:|\frac{1}{2\pi}\int_0^{2\pi}[h(e^{it})-h_0(e^{it})]g_j(e^{it})dt|<\epsilon,  j=1,2,\ldots, n, g_j\in L^1\right\}
\end{equation}
Let $\psi$ be a weak$^*$ generator of $H^{\infty}$. Then any set of the form \eqref{basis_topology} contains a polynomial in $\psi$ for any $h_0\in H^{\infty}$. Taking $g_j=\bar{f_j}$, $f_j\in H^1$, 
 we obtain that
 for any $h_0\in H^{\infty}(\D)$, for any $\epsilon>0$, and for any $f_j\in H^1(\D)$, $j=1,2,\ldots, n$, there is a polynomial $P(\psi)$ in $\psi$ such that 
 \begin{equation}\label{weak-star-neighborhood}
  |\frac{1}{2\pi}\int\limits_0^{2\pi}\left[P(\psi(e^{it}))-h_0(e^{it})\right]\overline{f_j(e^{it})}dt|<\epsilon, \quad j=1.2.\ldots, n.
 \end{equation}
Since $H^2(\D)\subset H^1(\D)$ we take $f_j\in H^2(\D)$ and \eqref{weak-star-neighborhood} gives 
\begin{equation*}
 |<P(\psi)-h_0, f_j>_{H^2}|<\epsilon,
\end{equation*}
that is the weak closure of polynomials in $\psi$ contains all $H^{\infty}$. In particular that closure contains all polynomials in $z$. Therefore weak closure of polynomials in $\psi$ coincides with $H^2(\D)$, hence the closure of polynomials in $\psi$ in $H^2(\D)$-norm is $H^2(\D)$. 
Consequently any weak$^*$ generator of $H^{\infty}(\D)$ is a generator of $H^2(\D)$. 

{\it 2.} Let $\vphi$ be a bounded univalent function and put $\Om=\vphi(\D)$. It is easily seen that the function $\vphi$ is a generator of $H^2(\D)$ if and only if the polynomials are dense in $H^2(\Om)$ in the norm of   $H^2(\Om)$. In \cite{Akr1} J.~Akeroyd constructed a bounded simply connected domain $\Om$ such that the polynomials are dense in $H^2(\Om)$. At the same time the mapping function $\vphi$ is not a weak$^*$ generator of $H^{\infty}(\D)$.  This statement follows from the fact that Akeroyd's domain is a subset of the unit disk $\D$ and the following result of D.~Sarason (see [8], Corollary 2, page 527):\\
{\it The function $\vphi$ fails to be a weak$^*$ generator if there is a domain $B$ containing $G$ properly such that
  \begin{equation*}
    \sup\limits_{z\in B}|f(z)|=\sup\limits_{z\in G}|f(z)|.
    \end{equation*}
  }
 \hfill$\Box$
\begin{remark}
Using similar arguments one can show that if $\psi$ is a weak$^*$ generator of $H^{\infty}$, then the  of the set of polynomials in $\psi$ is dense in $H^p$ for any $1<p<\infty$.
\end{remark}
\vskip 0,5truecm 
In the case of the Bergman space we note that the space $(A^1(\D))^*$ properly contains $H^{\infty}(\D)$. Let $\t$ be the topology on $H^{\infty}(\D)$ defined by the family of seminorms
\begin{equation}\label{Tau-topology}
 p_f(\psi)=|\int_{\D}\psi(z)\overline{f(\z)}d\s(z)|,\qquad \psi\in H^{\infty}(\D),\quad f\in A^1(\D),
\end{equation}
where $\s$ is the normalized Lebesque measure of the unit disk $\D$. A function $\psi\in H^{\infty}(\D)$ is called a $\t$-generator of $H^{\infty}(\D)$ if polynomials in $\psi$ are dense in $H^{\infty}$ in $\t-$topology.
Clearly  a $\t-$generator of $H^{\infty}$ is the Bergman space version of a weak$^*$ generator of $H^{\infty}(\D)$.

Applying the same arguments as above one concludes that every $\t-$ generator of $H^{\infty}$ is a generator of the Bergman space $A^2(\D)$ in the sense of Definition \ref{generator}. Suppose now that $\vphi$ is a generator of $A^2(\D)$ in the sense  
of Definition \ref{generator}. Therefore polynomials in $\vphi$ are dense in $A^2(\D)$ in the Bergman norm. In particular, polynomials in $\vphi$ are dense in $A^2(\D)$ in the weak topology. This means that for any $\epsilon>0$ and any $h_0\in H^{\infty}$, and for any finite set of functions $f_j\in A^2(\D)$ there is a polynomial $P$ such that
$$
|\int_{\D}(P(\vphi(z))-h(z))\overline{f_j(z)}d\s(z)|<\epsilon.
$$
Since $A^1(\D)$ contains functions that do not belong to $A^2(\D)$, 
 $\vphi$ is not a $\t-$generator of $H^{\infty}(\D)$.
 
 A simple example of a generator for the Hardy space $H^2(\D)$ and the Bergman space $A^2(\D)$ is a function 
 \begin{equation*}
 \vphi(z)=\frac{z_0-z}{1-\bar{z_0}z},\quad |z_0|<1.
 \end{equation*}
 The statement that $\vphi$ is a weak$^*$ generator of $H^{\infty}$ follows from the facts that $\vphi$ is a conformal mapping 
 of the unit disk onto itslef and that the unit disk is a Caratheodory domain. Direct proof of the fact that the linear hull of the functions $\vphi^k(z)$, $k=0,1,\ldots$ is dense in $H^2(\D)$ is also simple. Interesting examples of weak$^*$ generators, related to close-to-convex mapping are discussed in our article \cite{ABC}.
 
For the case of the Bergman space $A^2(\D)$ we need to show that from 
\begin{equation*}
\int\limits_{\D}f(z)\overline{\left(\frac{z_0-z}{1-\bar{z_0}z}\right)^k}dA(z)=0,\quad f\in A^2(\D),\quad k=0,1,2,\ldots,
\end{equation*}
it follows that $f=0$. After the change of variable $\z=\vphi(z)$ the last condition takes the form
\begin{equation*}
\int\limits_{\D}g(\z)\overline{\frac{\z^k}{(1-\bar{z_0}\z)^2}}dA(\z)=0, \quad k=0,1,\ldots,
\end{equation*}
where $g(\z)=-f(\vphi(\z))\vphi^{\prime}(\z)\in A^2(\D)$ and $\Vert f\Vert=\Vert g\Vert$, that is for any polynomial 
$P(\z)$
\begin{equation*}
\int\limits_{\D}g(\z)\overline{\frac{P(\z)}{(1-\bar{z_0}\z)^2}}dA(\z)=0, \quad k=0,1,\ldots.
\end{equation*}
Taking $P(\z)=\z^k(1-\bar{z_0}\z)^2$ one concludes that $g$ is orthogonal to $\z^k$ for $k=0,1,\ldots$. Hence $g=0$ and, consequently, $f=0$.
\vskip 0.5truecm 
A (closed) subspace of $\gH$ is called $\vphi$-invariant if it is invariant under the operator $M_{\vphi}$. We also denote by $Lat_{\gH}(\vphi)$ the lattice of $\vphi$-invariant subspaces of $\gH$.

For a set $S\subset \gH$ we denote $\left[S\right]_{\vphi}$ the smallest $\vphi-$ invariant subspace containing $S$. Similarly, $\left[S\right]_z$ is the smallest $z-$invariant subspace containing $S$. 

The next theorem gives a characterization of a generator in terms of its invariant subspaces.
\begin{theorem}\label{phi-invariance} Let $\vphi$ be a bounded univalent function in the unit disk $\D$. In order for $\vphi$ to be a generator it is necessary that every $\vphi-$invariant subspace $\gM\subset\gH$  such that $\ol{\gM\cap H^{\infty}}=\gM$ is also $z-$invariant, and sufficient that the invariant subspace 
$[\mathbf 1]_{\vphi}$ be $z-$invariant.

\end{theorem}
{\bf Proof.} At first we recall that a Toeplitz operator $T_g$ with symbol $g$ on $H^2(\D)$ ($A^2(\D)$) is defined by the formula 
$$
T_gh=P(gh)
$$ 
where $P$ is the orthogonal projection from $L^2(\T)$ onto $H^2(\D)$ (from $L^2(\D, dA)$ onto $A^2(\D)$). If $g\in H^{\infty}(\D)$ then the operator $T_g$ is bounded. 

Suppose $\vphi$ is a generator,  $f\in\gM\cap H^{\infty}$ and $h\in\gM^{\perp}$. Since $\gM$ is $\vphi-$invariant, $P(\vphi)f\in\gM$ for any polynomial  
$P$, that is 
$$
<P(\vphi)f,h>=0.
$$
The left side of the last equality can be written as $<P(\vphi),\bar{f}h>=<P(\vphi), T_{\bar{f}}h>$, where $T_{\bar{f}}$ is a  Toeplitz operator on $\gH$. Since $f\in H^{\infty}$ the operator $T_{\bar{f}}$ is bounded and $T_{\bar{f}}h\in\gH$. Thus we have 
$$
0=<P(\vphi)f,h>=<P(\vphi), T_{\bar{f}}h>.
$$
Now pick $\epsilon>0$. Because $\vphi$ is a generator there is a polynomial $P$ such that
$$
|<z-P(\vphi),T_{\bar{f}}h >|<\epsilon.
$$
Since $\epsilon>0$ is arbitrary it follows that
$$
0=<z,T_{\bar{f}}h>=<M_zf,h>.
$$
Hence $M_zf\in\gM$ for any $f\in\gM\cap H^{\infty}$ and, because  this intersection is dense in $\gM$ one concludes 
that $M_z\gM\subset\gM$. 

To prove the converse statement consider the $\vphi-$invariant subspace $[\mathbf 1]_{\vphi}$. Since this subspace is the closure in $\gH$ of polynomials in $\vphi$, functions from $H^{\infty}$ are dense in it. Consequently, $[\mathbf 1]_{\vphi}$ is also $z-$invariant. But any $z-$invariant subspace of $\gH$ which contains the function ${\mathbf 1}$ coincides with the whole $\gH$. Therefore $[\mathbf 1]_{\vphi}=\gH$, that is $\vphi$ is a generator.\hfill$\Box$
\begin{remark}
 It was proved earlier that a bounded univalent function $\psi$ is a weak$^*$ generator of $H^{\infty}$ if and only if 
 $Lat_{\gH}(\psi)=Lat_{\gH}(z)$. For the case $\gH=H^2(\D)$ it was proved by D.~Sarason in \cite{Sarason_1}, for 
 $\gH=A^2(\D)$   the statement was proved by P.~Bourdon in 
 \cite{Bourdon_1}.
\end{remark}


\begin{corollary} Let $\vphi$ be a generator of $\gH$ and let $\gM$ be a $\vphi-$invariant subspace.  If ${\rm dim}\; \gM^{\perp}<\infty$ then the subspace $\gM$ is $z$-invariant.
 \end{corollary}
The corollary follows from the fact that $H^{\infty}$ is dense in $\gH$ and the statement below that was proved for the more general situation in \cite{GK1}, Lemma 2.1:\\
{\it Let a Banach space $X$ be decomposed as the direct sum of a subspace $Y$ and a finite-dimensional subspace $Z$:
\begin{equation*}
X=Y\dot{+}Z,
\end{equation*}
 and $L$ is a dense linear subset of $X$. Then $Y\cap L$ is dense in $Y$.}\\

The following statement was proved in \cite{Carswell}, Lemma 4.1.
\begin{lemma}\label{Carswell_Lemma}
 Let $\vphi$ be a bounded univalent function on the unit disk $\D$ with $\vphi(0)=0$. If $\gM$ is a $z-$invariant subspace of $\gH$, then $\gM\ominus M_z\gM=\gM\ominus M_{\vphi}\gM$.
\end{lemma}
 
\begin{theorem} \label{wandering_property_th} Let $\vphi$ be a bounded univalent function on the unit disk $\D$ with $\vphi(0)=0$. Assume that $\vphi$ is a generator of $\gH$ and $\gM$ is a $\vphi-$ invariant subspace of $\gH$ such that:\\
(a) $\gM\cap H^{\infty}(\D)$ is dense in $\gM$;\\ 
(b) $(\gM\ominus\vphi\gM)\cap H^{\infty}(\D)$ is dense in $\gM\ominus\vphi\gM$.\\
Then 
 \begin{equation}\label{wandering_property}
  \gM=\left[\gM\ominus\vphi\gM\right]_{\vphi}.
 \end{equation}
\end{theorem}
{\bf Proof.} From Theorem \ref{phi-invariance} it follows that the subspace $\gM$ is $M_z-$invariant and Lemma \ref{Carswell_Lemma} gives that $\gM\ominus M_z\gM=\gM\ominus M_{\vphi}\gM$. If $\gH=H^2(\D)$ the fact that $\gM=\left[\gM\ominus M_z\gM\right]_z$ is the well known Wold decomposition. If $\gH=A^2(\D)$ then  
from the result of Aleman, Richter, and Sundberg(see \cite{ARS}) it follows that $\gM$ is the smallest $z-$invariant subspace that contains $\gM\ominus M_z\gM$. 

Observe that because $\vphi\in H^{\infty}(\D)$ we have 
$\left[\gM\ominus M_{\vphi}\gM\right]_{\vphi}\cap H^{\infty}(\D)$ is dense in $\left[\gM\ominus M_{\vphi}\gM\right]_{\vphi}$. Since the last subspace is $\vphi-$invariant we refer again to Theorem \ref{phi-invariance} and conclude that it is $z-$invariant, and is a subspace of $\gM$. Therefore $\left[\gM\ominus M_{\vphi}\gM\right]_{\vphi}=\gM$. \hfill$\Box$

\begin{remark} If in Theorem \ref{wandering_property_th} 
 $\dim(\gM\ominus\vphi\gM)<\infty$, then condition (a) implies condition (b).
\end{remark}

Let $\gM$ be a $\vphi-$invariant subspace of $\gH$. In what follows the quantity $\dim(\gM\ominus\vphi\gM)$ is called the $\vphi-$index of $\gM$. The $z-$index of a $z-$invariant subspace is defined similarly.

\begin{theorem}\label{1-generated_subspace} Let $\vphi$ be a generator of $\gH$ and $\vphi(0)=0$.  Assume $f\in(\gH\cap H^{\infty})$ and $\gM=[f]_{\vphi}$. Then the $\vphi-$index of $\gM$ equals 1 and $\gM=[f]_z$. Consequently, the $z-$index of $\gM$ is also 1.
\end{theorem}
{\bf Proof}. The proof of the first statement of the theorem is, in fact, repetition of the corresponding proof of Theorem 4, Chapter 8, of \cite{DS}. We include it for completeness.

Represent $f$ in the form $f=f_1+f_2$, where $f_1\in\gM\ominus\vphi\gM$ and $f_2\in\vphi\gM$. Choose an arbitrary $g\in\gM\ominus\vphi\gM$. Then $g\in[f]_{\vphi}$, consequently there is a sequence of polynomials of $\vphi$, say $Q_n(\vphi)$, such that $\Vert Q_n(\vphi)f-g\Vert\to 0$.
Define $h_n^{(1)}$ and and $h_n^{(2)}$ by the formulas
\begin{gather*}
 h_n^{(1)}=Q_n(\vphi(0))f_1-g\in\gM\ominus\vphi\gM,\\
 h_n^{(2)}=(Q_n(\vphi)-Q_n(\vphi(0)))f_1+Q_n(\vphi)f_2\in\vphi\gM.
\end{gather*}
The second inclusion follows from the fact that $Q_n(\vphi)-Q_n(\vphi(0))=(\vphi-\vphi(0))R_n(\vphi)$, $R_n(\phi)$ is a polynomial of $\vphi$ and $\vphi(0)=0$.
We have $h_n^{(1)}+h_n^{(2)}=Q_n(\vphi)f-g$, therefore
\begin{equation*}
 \Vert h_n^{(1)}\Vert^2+\Vert h_n^{(2)}\Vert^2=\Vert  Q_n(\vphi)f-g\Vert^2\to 0.
\end{equation*}
In particular, the property $\Vert h_n^{(1)}\Vert=\Vert Q_n(\vphi(0))f_1-g\Vert\to 0$ implies $g=\l f_1$, and therefore  
 the $\vphi-$index of $\gM$ is 1.

From our assumptions and Theorem \ref{phi-invariance} it follows that $\gM$ is $z-$invariant. Put $\gM^{\prime}=[f]_z$. Then $\gM^{\prime}\subset\gM$ and $\gM^{\prime}$ is $z-$invariant. Hence $\gM^{\prime}$ is also a $\vphi-$invariant subspace which contains $f$. But $\gM=[f]_{\vphi}$ is  {\it the smallest} $\vphi-$invariant subspace which contains $f$. Consequently, $\gM\subset\gM^{\prime}$. Hence $\gM=\gM^{\prime}$. 

Finally, any single-generated $z$-invariant subspace has $z-$index 1. \hfill$\Box$
\vskip 0.5truecm
\section{Density of $H^{\infty}(\D)$ in the Invariant Subspaces of a Generator}
The previous considerations raise the following questions:\\
{\it Let $\vphi$ be a generator of $\gH$ and $\gM$ is a $\vphi-$invariant subspace of $\gH$ such that {\rm dim}$(\gH\ominus\gM)=\infty$.\\ Is it true that
  $\gM\cap H^{\infty}\ne\{0\}$?\\
Is $\gM\cap H^{\infty}$ dense in $\gM$?}\\

For the case of $\gH=A^2(\D)$ the answer is given by the following theorem.
\begin{theorem}
 For any function $\psi$ analytic and bounded in the unit disk $D$ there is a singly-generated $\psi-$invariant subspace $\gM\subset A^2(\D)$ such that $\gM\cap H^{\infty}(\D)=\{0\}$. 
\end{theorem}
{\bf Proof.} Let $f\in A^2(\D)$ be such that its zeros do not satisfy the Blaschke condition and put $\gM=[f]_{\psi}$. Clearly, zeros of any function which belong to $\gM$ also do not satisfy the Blaschke condition, consequently, 
$\gM\cap H^{\infty}(\D)=\{0\}$.\hfill$\Box$
\vskip 0.5truecm
 Now we consider the case $\gH=H^2(\D)$. Let $f\in H^2(\D)\setminus H^{\infty}(\D)$. Then a function $f_T\in H^{\infty}$ is called $f-$truncating, if $f_Tf\in H^{\infty}$.
There are many ways to construct an $f-$truncating function. One of them is the following. 
 Define  a real valued function $\om(e^{it})$ as follows:
\begin{equation*}
 \om(e^{it})=\begin{cases}
              1& \text{if}\; |f(e^{it})|\le 1\\
              1/|f(e^{it})| &\text{if}\; |f(e^{it})|\ge 1.
             \end{cases}
\end{equation*}
We have $|\om(e^{it})|\le 1$ and $\log{\om}\in L^1(\T)$. Therefore there is an outer function $f_T\in H^{\infty}\subset H^2$ such that 
$|f_T(e^{it})|=\om(e^{it})$. It is obvious that $ff_T\in H^{\infty}(\D)$ and $\sup\{|f(z)f_T(z)|:|z|<1\}\le 1$.

Assume now that $\vphi$ is a generator. Therefore, there is a sequence of polynomials of $\vphi$, say $\{P_n(\vphi)\}$, such that $\Vert P_n(\phi)-f_T\Vert_2\to 0$. Such a sequence 
$\{P_n(\vphi)\}$ we shall call an $f-$truncating sequence. 
Since
$$
\int\limits_{0}^{2\pi}|P_n(\vphi(e^{it}))-f_T(e^{it})|^2dt\to 0
$$
there is a subsequence $P_{n_k}(\vphi)$ such that 
$P_{n_k}(\vphi(e^{it}))\to f_T(e^{it})$ for a.e. $t\in[0,2\pi]$.
Without loss of generality we may assume that $\{P_n(\vphi)\}$ converges to $f_T$ almost everywhere on $\T$. In particular $P_n(\vphi(e^{it}))f(e^{it})\to f_T(e^{it})f(e^{it})$ on the set of full measure on $\T$.
\begin{theorem}\label{contains_H_infty}
 Let $\vphi$ be a generator and  let $\gM$ be a $\vphi-$invariant subspace of the Hardy space $H^2(\D)$. 
Suppose $f\in\gM$ is not in $H^{\infty}(\D)$ and an $f$-truncating sequence $\{P_n(\vphi)\}$ is uniformly bounded in $H^{\infty}$-norm, that is $\sup_n{\Vert P_n(\vphi)\Vert_{\infty}}\le C<\infty$. Then $\gM\cap H^{\infty}(\D)\neq\{0\}$. 
 \end{theorem}
 {\bf Proof.} Let $f_T\in H^2\cap H^{\infty}$ be the $f-$ truncating function constructed in the previous paragraph. We need to show that $P_n(\vphi)f\to f_Tf$ in $H^2$-norm. Since 
 $P_n(\vphi)f\in\gM$ this will prove the statement.
 We have 
 \begin{gather*}
  \int_0^{2\pi}|P_n(\vphi(e^{it}))f(e^{it})-f_T(e^{it})f(e^{it})|^2dt=\\
\int_0^{2\pi}|P_n(\vphi(e^{it}))-f_T(e^{it})|^2|f(e^{it})|^2
dt.
\end{gather*}
The expression under integral sign converges to zero almost everywhere and is dominated by $C|f|^2$ where $C$ is a positive constant. Now Lebesgue's dominated convergence theorem gives the desired result.  \hfill$\Box$

\begin{theorem}\label{H-infty dense}
 Let $\psi$ be a weak-star generator of $H^{\infty}$ and let $\gM$ be a $\psi-$invariant subspace of $H^2(\D)$. Then $\gM\cap H^{\infty}$ is dense in $\gM$ in $H^2$ norm. 
\end{theorem}
{\bf Proof.} Let the functions $f\in\gM$ and $f_T$ be as in the previous theorem. For any $h\in H^2(\D)$ and any polynomial $P$ we have 
\begin{equation*}
 <P(\psi)f-f_Tf,h>=\frac{1}{2\pi}\int_0^{2\pi}(P(\psi(e^{it}))-f_T(e^{it}))f(e^{it})\overline{h(e^{it})}dt.
\end{equation*}
Now pick $\epsilon>0$. Since $f\bar{h}\in L^1(\T)$  and $\psi$ is a weak-star generator of $H^{\infty}$, there is a polynomial $P$ in $\psi$ such that the absolute value of the last integral is less than $\epsilon$. Since $P(\psi)f\in\gM$ it means that $f_Tf$ belongs to the weak closure of $\gM$. Since weak closure of a subspace coincides with its norm closure, it proves that $\gM\cap H^{\infty}\ne\{0\}$. 
 Put $\gM^{\prime}=c.l.h.\{P(\psi)f_Tf: 
  P\text{ is a polynomial}\}$.  Then the subspace
$\gM^{\prime}$ has the following properties:
\begin{enumerate}
 \item 
 $\gM^{\prime}\subset\gM;$
 \item $\gM^{\prime}$ is $\psi-$invariant;
 \item $\gM^{\prime}\cap H^{\infty}$ is dense in $\gM^{\prime}$ in $H^2$ norm.
\end{enumerate}
Let $\Om$ is the collection of all subspaces of $\gM$ which have these three properties. We have shown that $\Om\ne\emptyset$. Partially order $\Om$ by the set inclusion. By the Hausdorff's maximal principle there exists a maximal totally ordered subcollection $\Om^{\prime}$ of $\Om$. Denote by $\gM_0$ the union of all $\gM^{\prime}$, where $\gM^{\prime}$ is a member of $\Om^{\prime}$. Then $\gM_0$ is a maximal subspace of $\gM$ which satisfies the three conditions above.

We claim that $\gM_0=\gM$. Observe first that 
$\gM\setminus\gM_0$ does not contain any $H^{\infty}$ vectors. Suppose $\gM_0\ne\gM$ and pick $h\in\gM\ominus\gM_0$. Then $h\notin H^{\infty}$. Let $h_T$ be the $h$-truncating function. Then $h_Th\in H^{\infty}$, hence $h_Th\in\gM_0$ and $P(\psi)h_Th\in\gM_0$ for any polynomial $P$. Consequently, 
\begin{equation*}
 \frac{1}{2\pi}\int\limits_0^{2\pi}P(\psi(e^{it}))h_T(e^{it})|h(e^{it})|^2dt=0
\end{equation*}
for any polynomial $P$. Since $h_T|h|\in L^1(\T)$ and $\psi$ is a weak-star generator of $H^{\infty}$ from the last equality one easily deduces that 
\begin{equation*}
  \frac{1}{2\pi}\int\limits_0^{2\pi}e^{int}h_T(e^{it})|h(e^{it})|^2dt=0,\qquad n=0,1,2,\ldots
\end{equation*}\label{fmriesz}
from which it follows that
\begin{equation}
 h_T|h|^2=k,
\end{equation}
where $k\in H^1_0$. Therefore
\begin{equation*}
 |h|^2=\frac{k}{h_T}.
\end{equation*}
Since the left side is a real-valued function and right side is a function of the Nevanlina class and $k(0)=0$ one concludes that that equality \eqref{fmriesz} is possible only for $h=0$ and $k=0$. It proves our claim and the theorem. 
\hfill$\Box$ 

Combining the last theorem with Theorem \ref{phi-invariance} we obtain the following statement.
\begin{corollary}
Let $\psi$ be a bounded univalent function in the unit disk $\D$. Then $\psi$ is a weak$^*$ generator of $H^{\infty}$ if and only if every $\psi-$invariant subspace of $H^2(\D)$ is also $z-$invariant.
\end{corollary}
The statement of the Corollary was obtained by D. Sarason \cite{Sarason_1} using another method.




Valentin~V.~Andreev: Department of Mathematics, Lamar University, \\Beaumont, TX 77710, USA, \texttt{vvandreev@aol.com}\\
Miron~B.~Bekker: Department of Mathematics, The University of Pittsburgh at Johnstown, 450 Schoolhouse Rd, Johnstown, PA 15904, USA \texttt{bekker@pitt.edu}\\
Joseph~A.~Cima: Department of Mathematics, the University of \\North Carolina at Chapel Hill, CB 3250, 329 Phillips Hall, Chapel Hill, NC 27599, USA, \texttt{cima@email.unc.edu}
\end{document}